\documentclass[11pt,leqno]{article}

\usepackage{graphicx}
\usepackage{latexsym} 
\usepackage{amsmath}
\usepackage{amssymb}
\usepackage{amsfonts}
\usepackage{enumerate}
\usepackage{theorem}

\theoremstyle{plain}
\newtheorem{teo}{Theorem}[section]
\newtheorem{lem}[teo]{Lemma}
\newtheorem{cor}[teo]{Corollary}
\newtheorem{prop}[teo]{Proposition}

{\theorembodyfont{\rmfamily}
\newtheorem{defin}[teo]{Definition}
\newtheorem{oss}[teo]{Remark}

}

\renewcommand{\eqref}[1]{\textnormal{(\ref{#1})}}

\numberwithin{equation}{section}

\newcommand{\cvd}{\hfill$\square$}
\newcommand{\proof}[1]{\noindent\textsc{Proof#1}}
\newcommand{\rmi}{\mathrm{i}}
\newcommand{\rme}{\mathrm{e}}

\title{Stability for the acoustic scattering problem for sound-hard scatterers}

\author{Giorgio Menegatti\footnotemark[1]\hspace{1em}and\hspace{1em}Luca Rondi\footnotemark[2]}

\date{}

\begin{document}

\maketitle
\footnotetext[1]{Dipartimento di Matematica e Geoscienze, Universit\`a degli Studi di Trieste, via Valerio 12/1, 34127 Trieste ITALY and SISSA, via Bonomea 265, 34136 Trieste ITALY. E-mail: \texttt{giorgio.menegatti@teletu.it}}
\footnotetext[2]{Dipartimento di Matematica e Geoscienze, Universit\`a degli Studi di Trieste, via Valerio 12/1, 34127 Trieste ITALY. E-mail: \texttt{rondi@units.it}}

\setcounter{section}{0}
\setcounter{secnumdepth}{2}

\begin{abstract}
We study the stability for the direct acoustic scattering problem with sound-hard scatterers with minimal regularity assumptions on the scatterers. The main tool we use for this purpose is the convergence in the sense of Mosco.

We obtain uniform decay estimates for scattered fields and we investigate how a sound-hard screen may be approximated by thin sound-hard obstacles.

\medskip


\medskip


\end{abstract}

\section{Introduction}\label{intro}
We are interested in the following direct scattering problem. Let us send a time-harmonic acoustic  wave across a medium in $\mathbb{R}^N$, $N\geq 2$. If the medium is homogeneous such an incident wave remains unperturbed, whereas the presence of a scatterer creates a perturbation by producing a scattered wave. The scattered wave is radiating and it satisfies a suitable boundary condition on the boundary of the scatterer, that depends on the nature of the scatterer.

As incident wave we use a time-harmonic plane wave, which is characterized by
its field, referred to as the \emph{incident field}.
The incident field is given by $u^i(x)=\rme^{\rmi kx\cdot d}$, $x\in\mathbb{R}^N$, where
$k>0$ is the wavenumber and $d\in \mathbb{S}^{N-1}$ is the direction of
propagation.
A \emph{scatterer} $K$ in $\mathbb{R}^N$ is a compact subset of $\mathbb{R}^N$ such that $\mathbb{R}^N\backslash K$ is connected. We say that a scatterer is an \emph{obstacle} if
$K$ is the closure of its interior, whereas a scatterer is called a \emph{screen} if its interior is empty. If the incident wave encounters a scatterer $K$, then a scattered wave is created and its corresponding field will be called the \emph{scattered field} and denoted by $u^s$.
The \emph{total field} $u$ of the acoustic wave is given by the sum of the incident field
and the scattered field and solves the following direct scattering problem
\begin{equation}\label{dirpbm0}
\left\{\begin{array}{ll}
\Delta u + k^2u=0 & \text{in }\mathbb{R}^N\backslash K\\
u=u^i+u^s&\text{in }\mathbb{R}^N\backslash K\\
B.C.&\text{on }\partial K\\
\lim_{r\to \infty}r^{(N-1)/2}\left(\frac{\partial u^s}{\partial r}-\rmi ku^s\right)=0
& r=\|x\|,
\end{array}\right.
\end{equation}
where the last limit is the \emph{Sommerfeld radiation condition} and corresponds to the fact that
the scattered wave is radiating.
The boundary condition on the boundary of $K$ depends on the character  of the scatterer $K$. For instance, if $K$ is \emph{sound-soft}, then $u$ satisfies the following homogeneous Dirichlet condition
\begin{equation}\label{soundsoft}
u=0\quad\text{on }\partial K,
\end{equation}
whereas if $\Sigma$ is 
\emph{sound-hard} we have 
\begin{equation}\label{soundhard}
\nabla u\cdot \nu=0\quad\text{on }\partial K,
\end{equation}
that is a homogeneous Neumann condition.
Other conditions such as the impedance boundary condition or transmission conditions for penetrable scatterers may occur in the applications.

We recall that 
the Sommerfeld radiation condition implies that
the scattered field has the asymptotic behaviour of an outgoing spherical wave, namely
\begin{equation}\label{asympt0}
u^s(x)=\frac{\rme^{\rmi k\|x\|}}{\|x\|^{(N-1)/2}}\left\{u^s_{\infty}(\hat{x})
+O\left(\frac{1}{\|x\|}\right)\right\},
\end{equation}
where $\hat{x}=x/\|x\|\in \mathbb{S}^{N-1}$ and
$u^s_{\infty}$ is the so-called \emph{far-field pattern} of $u^s$. In particular, the scattered field satisfies the following decay property for some positive constants $E$ and $R$
\begin{equation}\label{decay0}
|u^s(x)|\leq E\|x\|^{-(N-1)/2}\quad\text{for any }x\in\mathbb{R}^N\text{ so that }\|x\|\geq R.
\end{equation}

We refer to \cite{Wil} for further details about the direct scattering problem
\eqref{dirpbm0}. See also \cite{Col e Kre98} where the corresponding inverse problems are considered.

Here we wish to find suitable conditions on sound-hard scatterers $K$ that make the scattering problem \eqref{dirpbm0} stable with respect to $K$, that is with respect to variations
in the scatterer. Such a problem is of interest 
since in general solutions to Neumann problems for elliptic equations may not be stable under domain variations as the so-called Neumann sieve case points out, see for instance \cite{Mur}.

We recall that the simpler case of sound-soft scatterers has been treated in \cite{Ron03}. There,  under minimal regularity assumptions on the scatterers, stability results for the solutions of the direct scattering problems with respect to the scatterer and uniform decay estimates for the corresponding scattered fields have been obtained.

In order to evaluate distances between scatterers and convergence of scatterers we shall use the Hausdorff distance. In particular, a sequence of bounded open sets $\Omega_n$ contained in $B_R$ for some $R>0$ converges to an open set $\Omega$  if their complements in $\overline{B_R}$ converges to the complement of $\Omega$ in the Hausdorff distance.
The main tool we shall use is the so-called Mosco convergence of the corresponding Sobolev spaces $H^1(\Omega_n)$ to $H^1(\Omega)$ . It has been proved that Mosco convergence is essentially equivalent to convergence of solutions of Neumann problems, at least for elliptic equations which are strictly coercive in $H^1$. We shall show that Mosco convergence is indeed a sufficient condition also for the convergence of solutions  of Neumann problems for the Helmholtz equation, and in particular for scattering problems, provided a uniform Sobolev type inequality holds true for any $\Omega_n$. These results are contained in Section~\ref{stabsec}.

Then we are interested in finding classes of scatterers for which we have uniform decay estimates for the scattered fields, that is estimates like \eqref{decay0} with constants $E$ and $R$ independent on $K$. This will be done in Section~\ref{scatsec}.
Again the key point, besides a uniform Sobolev type inequality, is establishing sufficient conditions on scatterers that guarantee, under convergence in the Hausdorff distance, convergence of corresponding Sobolev spaces in the sense of Mosco. Such a problem, motivated by convergence of solutions of Neumann problems, has been extensively studied in the literature.

In dimension $2$, the problem is fully solved since Bucur and Varchon gave a necessary and sufficient condition, \cite{Buc-Var}. The starting point of \cite{Buc-Var} was the sufficient condition proved by Chambolle and Doveri, \cite{Cha-Dov}, which is still a fairly convenient one to use for the applications. In dimension $2$ complex analytic techniques, in particular duality arguments, are crucial for obtaining these results.

In dimension $3$ and higher, there is a result by Giacomini, \cite{Gia}, where the admissible sets $K$ are Lipschitz hypersurfaces with Lipschitz boundaries. 

We precisely recall these assumptions in dimension $2$ and in dimension $3$ and higher in Section~\ref{scatsec}. Moreover we introduce a new sufficient condition, see Theorem~\ref{Moscoconvergencetheo}, which is based on a generalization of a class of sets previously
introduced in \cite{Ron06}. These new admissible sets still consists of Lipschitz hypersurfaces with Lipschitz boundaries, with a definition which is slightly less general than that given by Giacomini. The main advantage is that we allow the presence of more hypersurfaces and especially that they may intersect in a transversal way.

Finally, in Section~\ref{examplesec}, we show how to approximate a (Lipschitz) screen by a thin obstacle surrounding it. Such a result may ease the numerical computation of the solution to the scattering problem for a screen and it is also related to some issues arising in the so-called cloaking problem.

The plan of the paper is the following.
In Section~\ref{stabsec} we deal with the stability for the direct acoustic scattering problem with sound-hard scatterers. We first discuss Mosco convergence and the stability of Neumann problems for the Helmholtz equation, then we treat scattering problems. In Section~\ref{scatsec} we define suitable classes of sound-hard scatterers whose corresponding scattered fields satisfy a uniform decay estimate. Finally in Section~\ref{examplesec} we show how we can approximate sound-hard screens by thin sound-hard obstacles.

\subsubsection*{Acknowledgements} LR is partly supported by Universit\`a
degli Studi di Trieste through Finanziamento per Ricercatori di Ateneo 2009 and by GNAMPA, INdAM, through 2012 projects. 

\section{Stability of Neumann problems for the Helmholtz equation}\label{stabsec}

For any $x\in\mathbb{R}^N$, $N\geq 2$,  we denote $x=(x',x_N)\in\mathbb{R}^{N-1}\times \mathbb{R}$ and $x=(x'',x_{N-1},x_N)\in\mathbb{R}^{N-2}\times\mathbb{R}\times\mathbb{R}$.
For any $s>0$ and any $x\in\mathbb{R}^N$, 
$B_s(x)$ denotes the ball contained in $\mathbb{R}^N$ with radius $s$ and center $x$, whereas $B'_s(x')$ denotes the ball contained in $\mathbb{R}^{N-1}$ with radius $s$ and center $x'$.
Moreover, $B_s=B_s(0)$ and $B'_s=B'_s(0)$. Finally, for any $E\subset \mathbb{R}^N$, we denote $B_s(E)=\bigcup_{x\in E}B_s(x)$.

\subsection{Mosco convergence}\label{Moscosubsec}

We recall that, given $\{A_n\}_{n\in\mathbb{N}}$, a sequence of closed subspaces of a reflexive Banach space $X$,
we denote
$$A'=\{x\in X:\ x=w\text{-}\lim_{k\to\infty} x_{n_k},\ x_{n_k}\in A_{n_k}\}$$
and
$$A''=\{x\in X:\ x=s\text{-}\lim_{n\to\infty} x_n,\ x_n\in A_n\}.$$
We note that $A'$ and $A''$ are subspaces of $X$, that $A''\subset A'$, and that $A''$ is closed. We say that
$A_n$ converges, as $n\to\infty$, to a closed subspace $A$ in the sense of Mosco if $A=A'=A''$. In other words, the following two conditions need to be satisfied.
\begin{enumerate}[i)]
\item For any $x\in X$, if there exists a subsequence $A_{n_k}$ and a sequence $x_k$, $k\in\mathbb{N}$,
such that $x_k$ converges weakly to $x$ as $k\to\infty$ and $x_k\in A_{n_k}$ for any $k\in\mathbb{N}$, then $x\in A$.
\item For any $x\in A$, there exists a sequence $x_n\in A_n$, $n\in\mathbb{N}$, converging strongly to $x$
as $n\to\infty$.
\end{enumerate}

Let $\Omega_1\Subset\Omega$ be bounded open sets contained in $\mathbb{R}^N$, $N\geq 2$. We assume that $\Omega$ has a Lipschitz boundary. We denote
$$\mathcal{A}=\{K\subset \overline{\Omega_1}:\ K\text{ is compact}\}.$$
We have that $\mathcal{A}$ is compact with respect to the Hausdorff distance, see for instance \cite{DalMaso}.

Let us notice that if $K_n\in\mathcal{A}$, $n\in\mathbb{N}$, converges in the Hausdorff  distance, as $n\to\infty$, to $K\in\mathcal{A}$, then we also have $\lim_{n\to\infty}|K_n\backslash K|=0$. 

For any $K\in\mathcal{A}$, we consider the isometric immersion of $H^1(\Omega\backslash K)$ into $L^2(\Omega,\mathbb{R}^{N+1})$ as follows. To each $u\in H^1(\Omega\backslash K)$ we associate the vector $(u,\nabla u)$ with the convention that 
$u$ and $\nabla u$ are extended to zero in $K$. In such a way we may consider $H^1(\Omega\backslash K)$ as a closed subspace of $L^2(\Omega,\mathbb{R}^{N+1})$.

Given a sequence $\{K_n\}_{n\in\mathbb{N}}$ contained in $\mathcal{A}$ and $K\in\mathcal{A}$, we say that $H^1(\Omega\backslash K_n)$ converges, as $n\to\infty$, to $H^1(\Omega\backslash K)$ \emph{in the sense of Mosco}
if this holds considering $H^1(\Omega\backslash K_n)$, $n\in\mathbb{N}$, and $H^1(\Omega\backslash K)$ as subspaces of 
$L^2(\Omega,\mathbb{R}^{N+1})$.

We are interested to find sufficient conditions on $K_n$, $n\in\mathbb{N}$, and $K$  in order to have that 
$H^1(\Omega\backslash K_n)$ converges, as $n\to\infty$, to $H^1(\Omega\backslash K)$ in the sense of Mosco. The following results are present in the literature, under the assumption that $K_n$ converges, as $n\to\infty$, to $K$ in the Hausdorff distance. For any $N\geq 2$, a sufficient condition on $K_n$, $n\in\mathbb{N}$, is established in \cite{Gia}. For $N=2$, instead, the key assumption is that the number of connected components of $K_n$, $n\in\mathbb{N}$, is uniformly bounded. Under this assumption, and the convergence in the Hausdorff distance, Bucur and Varchon, \cite{Buc-Var,Buc-Var2}, proved that  the convergence in the sense of Mosco holds if and only if  $|\Omega\backslash K_n|$ converges, as $n\to\infty$, to $|\Omega\backslash K|$.
A sufficient condition, and a starting point for the result of Bucur and Varchon, has been 
given by Chambolle and Doveri, \cite{Cha-Dov}. They proved the convergence in the sense of Mosco provided we have a uniform bound on the number of connected components of $\partial K_n$ and a uniform bound on $\mathcal{H}^1(\partial K_n)$, $n\in\mathbb{N}$.

Let us remark that the convergence in the sense of Mosco is built in such a way that is essentially equivalent to convergence of solutions of Neumann problems in the following sense. For any $f\in L^2(\Omega)$ and for any $K\in \mathcal{A}$, let $u=u(K,f)$ be the solution to the following problem
\begin{equation}\label{ellipticequation}
\left\{
\begin{array}{ll}
-\Delta u+u=f & \text{in }\Omega\backslash K\\
\nabla u\cdot\nu=0 & \text{on }\partial (\Omega\backslash K),
\end{array}\right.
\end{equation}
that is $u\in H^1(\Omega\backslash K)$ and satisfies
$$\int_{\Omega\backslash K}\nabla u\cdot \nabla  \varphi+\int_{\Omega\backslash K}u\varphi=\int_{\Omega\backslash K}f\varphi\quad\text{for any }\varphi\in H^1(\Omega\backslash K).$$
Then the following result holds.

\begin{prop}\label{Moscocaratt}
Let us fix  a sequence $\{K_n\}_{n\in\mathbb{N}}$ contained in $\mathcal{A}$ and $K\in\mathcal{A}$.

If $H^1(\Omega\backslash K_n)$ converges, as $n\to\infty$, to $H^1(\Omega\backslash K)$ in the sense of Mosco, then for any $f\in L^2(\Omega)$ we have that
$u_n=u(K_n,f)$, solution to \eqref{ellipticequation} with $K$ replaced by $K_n$, converges to $u=u(K,f)$, solution to \eqref{ellipticequation}, where as before the convergence is in the sense of 
$L^2(\Omega,\mathbb{R}^{N+1})$.

Conversely, if $K_n$ converges, as $n\to\infty$, to $K$ in the Hausdorff distance and for any $f\in L^2(\Omega)$ we have that
$u_n=u(K_n,f)$ converges to $u=u(K,f)$ in $L^2(\Omega,\mathbb{R}^{N+1})$, then
$H^1(\Omega\backslash K_n)$ converges to $H^1(\Omega\backslash K)$ in the sense of Mosco.
\end{prop}

\proof{.} See for instance Proposition~3.2 and Remark 3.3 in \cite{Buc-Var2}.\cvd

\bigskip

By Proposition~\ref{Moscocaratt}, taking $f\equiv 1$, we easily infer that the convergence of $|\Omega\backslash K_n|$ to $|\Omega\backslash K|$, as $n\to\infty$, is a necessary condition for the convergence of  $H^1(\Omega\backslash K_n)$ to $H^1(\Omega\backslash K)$ in the sense of Mosco, for any $N\geq 2$. In fact,
under the convergence, as $n\to\infty$, in the sense of Mosco of $H^1(\Omega\backslash K_n)$ to $H^1(\Omega\backslash K)$, we have $\lim_{n\to\infty}|(\Omega\backslash K)\Delta(\Omega\backslash K_n)|=0$ or, equivalently,  $\lim_{n\to\infty}|K\Delta K_n|=0$, which implies that
$|\Omega\backslash K_n|$ converges, as $n\to\infty$, to $|\Omega\backslash K|$.

Moreover, the following important observation will be of use.

\begin{prop}\label{convboundary}
Let $\{K_n\}_{n\in\mathbb{N}}$ be a sequence in $\mathcal{A}$ and let $K$ and $\tilde{K}\in\mathcal{A}$. Let $\tilde{A}_n=H^1(\Omega\backslash \partial K_n)$, $n\in\mathbb{N}$, and $\tilde{A}=H^1(\Omega\backslash \tilde{K})$.
Let us assume that, as $n\to\infty$,  $\tilde{A}_n$ converges to $\tilde{A}$ in the sense of Mosco, 
that $\partial K_n$ converges to $\tilde{K}$ in the Hausdorff distance, and
that $K_n$ converges to $K$ in the Hausdorff distance.

Then $H^1(\Omega\backslash K_n)$ converges, as
$n\to \infty$, to $H^1(\Omega\backslash K)$ in the sense of Mosco.
\end{prop}

\proof{.} For any fixed $f\in L^2(\Omega)$, let 
$u_n=u(K_n,f)$, $n\in\mathbb{N}$, solution to \eqref{ellipticequation} with $K$ replaced by $K_n$, and $u=u(K,f)$, solution to \eqref{ellipticequation}.

For any $n\in\mathbb{N}$, since $u_n$ is zero inside $K_n$, we have that $u_n\in H^1(\Omega\backslash \partial K_n)$ and $u_n$ solves \eqref{ellipticequation} also with
$K$ replaced by $\partial K_n$.

By the first part of Proposition~\ref{Moscocaratt} we have that $u_n$ converges in 
$L^2(\Omega,\mathbb{R}^{N+1})$
to a function $\tilde{u}$ solving \eqref{ellipticequation} with $K$ replaced by $\tilde{K}$. We need to show that $\tilde{u}$ solves \eqref{ellipticequation} also for $K$, that is $\tilde{u}=u$. Then the proof is concluded by using the second part of Proposition~\ref{Moscocaratt}.

We observe that $\partial K\subset \tilde{K}\subset K$. First of all, we show that $\tilde{u}\in H^1(\Omega\backslash K)$, that is we need to show that $\tilde{u}$, and its gradient, are zero almost everywhere in $K$. It is sufficient to prove that for $K\backslash \tilde{K}=\stackrel{\circ}{K}\backslash \tilde{K}$, actually for any ball $B_\delta(x)$, with $\delta>0$ and such that $B_{2\delta}(x)\subset \stackrel{\circ}{K}\backslash \tilde{K}$. By the Hausdorff convergence, we have that the intersection of $B_\delta(x)$ with $\partial K_n$ is empty, for any $n$ large enough. Since $x\in K$, we deduce that  $B_\delta(x)$ is contained in $K_n$, again for any $n$ large enough. Finally, since $u_n$ is zero in $K_n$, we conclude that $\tilde{u}$ is zero in  $B_\delta(x)$. An analogous reasoning holds for the gradient of $\tilde{u}$, hence $\tilde{u}\in H^1(\Omega\backslash K)$.

Then, take $\varphi\in H^1(\Omega\backslash K)$. We prove that $\varphi$ belongs to $H^1(\Omega\backslash \tilde{K})$. Since $\tilde{K}\subset K$, we have that $\varphi$ and its gradient are zero in $\tilde{K}$. We notice that $\Omega\backslash \tilde{K}=(\Omega\backslash K)\cup(\stackrel{\circ}{K}\backslash \tilde{K})$. Since $\varphi\in H^1(\Omega\backslash K)$ and $\varphi$ and its gradient are zero in $\stackrel{\circ}{K}\backslash \tilde{K}$, we infer that $\varphi\in H^1(\Omega\backslash \tilde{K})$. So we notice that
$$\int_{\Omega\backslash K}f\varphi=\int_{\Omega\backslash \tilde{K}}f\varphi=
\int_{\Omega\backslash \tilde{K}}\nabla \tilde{u}\cdot \nabla  \varphi+\int_{\Omega\backslash \tilde{K}}\tilde{u}\varphi=\int_{\Omega\backslash K}\nabla \tilde{u}\cdot \nabla  \varphi+\int_{\Omega\backslash K}\tilde{u}\varphi$$
and the proof is concluded.\cvd

\bigskip

We conclude this subsection by pointing out the following two lemmas.

\begin{lem}\label{primolemma}
Let $\{K_n\}_{n\in\mathbb{N}}$ be a sequence in $\mathcal{A}$ converging, as $n\to\infty$, to $K\in\mathcal{A}$ in the Hausdorff distance. Let $A_n=H^1(\Omega\backslash K_n)$, $n\in\mathbb{N}$, and $A=H^1(\Omega\backslash K)$.

If $|K\backslash K_n|$ goes to zero as $n\to\infty$, then
we have that
$A'\subset A$.
\end{lem}

\proof{.} Let $\{A_{n_K}\}_{k\in\mathbb{N}}$ be a subsequence and let $(u_k,\nabla u_k)\in A_{n_k}$ for any $k\in\mathbb{N}$. We assume that $(u_k,\nabla u_k)$ weakly converges, as $k\to\infty$,
to $(v,V)$.

For any $x\in \Omega\backslash K$ there exist $r>0$ and $\overline{k}\in\mathbb{N}$ such that
$\overline{B_r(x)}\subset \Omega\backslash K_{n_k}$ for any $k\geq \overline{k}$. It follows immediately that $v\in H^1(B_r(x))$ and that $V=\nabla v$ in $B_r(x)$. Then it is easy to conclude that $v\in H^1(\Omega\backslash K)$ and $V=\nabla v$ in  $\Omega\backslash K$.

It remains to prove that $v$ and $V$ are $0$ almost everywhere in $K$. Clearly $(u_k,\nabla u_k)$ are identically equal to zero on $K\cap K_{n_k}$, for any $k\in \mathbb{N}$. Therefore, for any $\varphi\in L^2(K)$ we have
$$\left|\int_K u_k\varphi\right|= \left|\int_{K\backslash K_{n_k}} u_k\varphi\right|\leq \|u_k\|_{L^2(\Omega)}\left(\int_{K\backslash K_{n_k}}\varphi^2\right)^{1/2},$$
and clearly the right-hand side goes to zero as $k\to\infty$ since $|K\backslash K_{n_k}|$ goes to zero. Therefore $u_k$ weakly converges in $L^2(K)$, as $k\to\infty$, to zero, that is $v=0$ almost everywhere in $K$. A similar reasoning holds for $\nabla u_k$,  so
the proof is concluded.\cvd

\begin{lem}\label{secondolemma}
Let $\{K_n\}_{n\in\mathbb{N}}$ be a sequence in $\mathcal{A}$ and let $K\in\mathcal{A}$. Let $A_n=H^1(\Omega\backslash K_n)$, $n\in\mathbb{N}$, and $A=H^1(\Omega\backslash K)$.

If $K\subset K_n$ for any $n\in\mathbb{N}$ and $|K_n\backslash K|$ goes to zero as $n\to\infty$, then we have that $A\subset A''$. 
\end{lem}

\proof{.} Fixed $u\in H^1(\Omega\backslash K)$, just take, for any $n\in\mathbb{N}$, $u_n=u (1-\chi_{K_n})$. Clearly $u_n$ belongs to $H^1(\Omega\backslash K_n)$ and $(u_n,\nabla u_n)$ strongly converges to $(u,\nabla u)$ in $L^2(\Omega,\mathbb{R}^{N+1})$.\cvd

\begin{oss}
We note that, in the previous lemma, it would be sufficient to assume for instance that $K_n$ converges, as $n\to\infty$, to $K$ in the Hausdorff distance. In fact convergence in the Hausdorff distance implies that $|K_n\backslash K|$ goes to zero as $n\to\infty$.
\end{oss}

\begin{cor}
Let $\{K_n\}_{n\in\mathbb{N}}$ be a sequence in $\mathcal{A}$ converging, as $n\to\infty$, to $K\in\mathcal{A}$ in the Hausdorff distance. Let $A_n=H^1(\Omega\backslash K_n)$, $n\in\mathbb{N}$, and $A=H^1(\Omega\backslash K)$.

If $K\subset K_n$ for any $n\in\mathbb{N}$, then we have that $A_n$ converges to $A$, as $n\to\infty$, in the sense of Mosco.
\end{cor}

\subsection{Application to stability of Neumann problems}\label{stabsubsec}

We observe that Mosco convergence is enough to pass to the limit in Neumann problems, at least in a weak sense and under a mild compactness assumption. Let us fix $k>0$.

\begin{prop}\label{firstprop}
Let $\{K_n\}_{n\in\mathbb{N}}$ be a sequence in $\mathcal{A}$ and let $K\in\mathcal{A}$. Let $A_n=H^1(\Omega\backslash K_n)$, $n\in\mathbb{N}$, and $A=H^1(\Omega\backslash K)$.
Let us assume that, as $n\to\infty$,  $A_n$ converges to $A$ in the sense of Mosco.

For any $n\in\mathbb{N}$, let $u_n\in H^1(\Omega\backslash K_n)$ solve, in a weak sense,
\begin{equation}\label{u_neq}
\left\{\begin{array}{ll}
\Delta u_n+k^2u_n=0 &\text{in }\Omega\backslash K_n\\
\nabla u_n\cdot \nu=0 &\text{on }\partial K_n,
\end{array}\right.
\end{equation}
that is 
$$\int_{\Omega\backslash K_n}\nabla u_n\cdot\nabla \varphi-k^2 u_n\varphi=0
$$
for any $\varphi\in H^1(\Omega\backslash K_n)$ such that the support of $\varphi$ is compactly contained in $\Omega$.

Let us assume that, for a constant $C$,
\begin{equation}\label{boundedness}
\|u_n\|_{L^2(\Omega\backslash K_n)}\leq C\quad\text{for any }n\in \mathbb{N}.
\end{equation}

Then, up to a subsequence, we have that $u_n$ converges weakly in $L^2(\Omega)$ to
a function $u$ solving in the same weak sense
\begin{equation}\label{ueq}
\left\{\begin{array}{ll}
\Delta u+k^2u=0 &\text{in }\Omega\backslash K\\
\nabla u\cdot \nu=0 &\text{on }\partial K.
\end{array}\right.
\end{equation}
\end{prop}

\proof{.} Let us assume that, up to a subsequence, we have that $u_n$ is converging weakly to $u$
in $L^2(\Omega)$ as $n\to\infty$.

Let us fix an open subset $D$ such that $D$ is of class $C^1$, and  $\Omega_1\Subset D\Subset\Omega$. Let $r>0$ be such that $B_r(\partial D)\Subset \Omega\backslash\overline{\Omega_1}$.

By standard regularity estimates, we may infer that there exists a constant $C_1$ such that
$$\|u_n\|_{C^1(B_{r/2}(\partial D))}\leq C_1\quad\text{for any }n\in\mathbb{N}.$$
Therefore, since
$$\int_{D\backslash K_n}|\nabla u_n|^2=k^2\int_{D\backslash K_n} u_n^2+\int_{\partial D}(\nabla u_n\cdot\nu) u_n,$$
we conclude that
$\{(u_n,\nabla u_n)\}_{n\in\mathbb{N}}$ is uniformly bounded in $L^2(D,\mathbb{R}^{N+1})$, so we may assume that, up to a subsequence, it converges weakly to $(v,V)$ as $n\to\infty$. Clearly $v=u$, and,
by the Mosco convergence, we have that $u\in H^1(D\backslash K)$ and $V=\nabla u$.

Then, take $\varphi\in H^1(\Omega\backslash K)$ such that $\varphi=0$ outside $D$. By the Mosco convergence, we can find $\varphi_n\in H^1(\Omega\backslash K_n)$, such that
$\varphi_n=0$ outside $D$ for any $n\in\mathbb{N}$, that
 converges strongly to $\varphi$ as $n\to\infty$. Then, since
$$\int_{\Omega}\nabla u_n\cdot\nabla \varphi_n-k^2 u_n\varphi_n=0,
$$
we are able to pass to the limit and prove that 
$$\int_{\Omega}\nabla u\cdot\nabla \varphi-k^2 u\varphi=0.
$$
The proof immediately follows by changing in a suitable way the set $D$.\cvd

\bigskip

By Proposition~\ref{firstprop} and Proposition~\ref{convboundary}, we immediately infer this useful modification.

\begin{prop}\label{firstpropmodified}
Let $\{K_n\}_{n\in\mathbb{N}}$ be a sequence in $\mathcal{A}$ and let $K$ and $\tilde{K}\in\mathcal{A}$. Let $\tilde{A}_n=H^1(\Omega\backslash\partial K_n)$, $n\in\mathbb{N}$, and $\tilde{A}=H^1(\Omega\backslash \tilde{K})$.
Let us assume that, as $n\to\infty$,  $\tilde{A}_n$ converges to $\tilde{A}$ in the sense of Mosco, 
that $\partial K_n$ converges to $\tilde{K}$ in the Hausdorff distance, and
that $K_n$ converges to $K$ in the Hausdorff distance.

For any $n\in\mathbb{N}$, let $u_n\in H^1(\Omega\backslash K_n)$ solve \eqref{u_neq} and
let us assume that, for a constant $C$, \eqref{boundedness} holds.

Then, up to a subsequence, we have that $u_n$ converges weakly in $L^2(\Omega)$ to
a function $u$ solving \eqref{ueq}.
\end{prop}

In certain cases, weak convergence might not be enough. In order to have strong convergence in $L^2$, we need to modify the previous propositions in the following way.

\begin{prop}\label{secondprop}
Let $\{K_n\}_{n\in\mathbb{N}}$ be a sequence in $\mathcal{A}$ and let $K\in\mathcal{A}$.
Let $A=H^1(\Omega\backslash K)$, and, for any $n\in\mathbb{N}$,
$A_n=H^1(\Omega\backslash K_n)$ and $\tilde{A}_n=H^1(\Omega\backslash \partial K_n)$. Let 
us assume that, as $n\to\infty$,
$K_n$ converges to $K$ in the Hausdorff distance.

Let us assume that, as $n\to\infty$, either $A_n$ converges to $A$ in the sense of Mosco,
or $\tilde{A}_n$ converges to $\tilde{A}=H^1(\Omega\backslash \tilde{K})$ in the sense of Mosco
and $\partial K_n$ converges to $\tilde{K}$ in the Hausdorff distance, for some $\tilde{K} \in\mathcal{A}$. 


For any $n\in\mathbb{N}$, let $u_n\in H^1(\Omega\backslash K_n)$ solve \eqref{u_neq} and
assume that, for a constant $C$,
\begin{equation}\label{boundedness2}
\|(u_n,\nabla u_n)\|_{L^2(\Omega,\mathbb{R}^{N+1})}\leq C\quad\text{for any }n\in \mathbb{N}.
\end{equation}

Let us further assume that there exist constants $p>2$ and $C_1>0$ such that for any $n\in\mathbb{N}$ we have
\begin{equation}\label{Sobolev}
\|u\|_{L^p(\Omega\backslash K_n)}\leq C_1\|u\|_{H^1(\Omega\backslash K_n)}\quad\text{for any }u\in H^1(\Omega\backslash K_n).
\end{equation}

Then, up to a subsequence, we have that $u_n$ converges strongly in $L^2(\Omega)$ to
a function $u$ solving \eqref{ueq}.
\end{prop}

\proof{.} Either by Proposition~\ref{firstprop} or  by Proposition~\ref{firstpropmodified}, up to a subsequence, we have that $u_n$ converges
weakly in $L^2(\Omega)$ to
a function $u$ solving \eqref{ueq}.
By standard regularity estimates, and by the convergence in the Hausdorff distance, we may also assume that $u_n$ converges to $u$ strongly in $L^2$ on any compact subset of $\Omega\backslash K$.

For any $\varepsilon>0$, we can find $D_{\varepsilon}$, a compact subset of $\Omega\backslash K$, such that, denoting $E_{\varepsilon}=(\Omega\backslash K)\backslash D_{\varepsilon}$, we have
$|E_{\varepsilon}|\leq \varepsilon$.

Let us notice that, for any $n\in\mathbb{N}$ and any $\varepsilon>0$, fixed $p>2$, we have
$$\|u_n\|_{L^2(E_{\varepsilon})}\leq \|u_n\|_{L^p(E_{\varepsilon}\cap (\Omega\backslash K_n))}\varepsilon^{(p-2)/(2p)}.$$
Then, fixed $\varepsilon>0$ and $n\in\mathbb{N}$ such that $K_n\subset K\cup E_{\varepsilon}$,
we have that
 $$\|u_n-u\|_{L^2(\Omega)}\leq \|u\|_{L^2(E_{\varepsilon})}+\|u_n\|_{L^2(E_{\varepsilon})}+
\|u_n\|_{L^2(K\backslash K_n)}+ 
  \|u_n-u\|_{L^2(D_{\varepsilon})}.$$
By \eqref{Sobolev} and \eqref{boundedness2}, we infer that
$$\|u_n-u\|_{L^2(\Omega)}\leq \|u\|_{L^2(E_{\varepsilon})}+CC_1(\varepsilon^{(p-2)/(2p)}+|K\backslash K_n|^{(p-2)/(2p)})+ 
  \|u_n-u\|_{L^2(D_{\varepsilon})}.$$
Fixed $\delta>0$, we find $\varepsilon>0$ and $\overline{n}\in\mathbb{N}$ such that for any $n\geq \overline{n}$ we have  $K_n\subset K\cup E_{\varepsilon}$ and
$$\|u\|_{L^2(E_{\varepsilon})}+CC_1(\varepsilon^{(p-2)/(2p)}+|K\backslash K_n|^{(p-2)/(2p)})\leq \delta/2.$$
Here we have used the fact that Mosco convergence implies that $|K\Delta K_n|$ goes to $0$ as $n\to\infty$. 
Then by the convergence on compact subsets, there exists $\tilde{n}\geq \overline{n}$ such that for every $n\geq\tilde{n}$ we have  $\|u_n-u\|_{L^2(D_{\varepsilon})}\leq \delta/2$, therefore for every $n\geq\tilde{n}$ we have
$$\|u_n-u\|_{L^2(\Omega)}\leq \delta$$
and the proof is concluded.\cvd

\begin{oss} We observe that all the results in this subsection remain valid if we let the wavenumber $k$ depend on $n$, with the assumption that $k_n\geq 0$ converges to a real number $k_{\infty}\geq 0$. 
\end{oss}

We conclude this subsection finding sufficient conditions for the assumption \eqref{Sobolev} to hold. Given $\mathcal{C}$ a  fixed cone in $\mathbb{R}^N$, we say that an open set $D\subset\mathbb{R}^N$ satisfies the \emph{cone condition with cone} $\mathcal{C}$ if for
every $x\in D$ there exists a cone $\mathcal{C}(x)$ with vertex in $x$ and congruent to $\mathcal{C}$ such that $\mathcal{C}(x)\subset D$. We remark that by a cone we always mean a bounded not empty open cone.
We shall use two different kind of conditions. 

\begin{lem}\label{conelemma}
Let $\mathcal{C}$ be a fixed cone in $\mathbb{R}^N$.
Let $D$ be a bounded open set satisfying the cone condition with cone $\mathcal{C}$.

Then there exist constants $p>2$ and $C_1>0$ such that
\begin{equation}\label{SobolevD}
\|u\|_{L^p(D)}\leq C_1\|u\|_{H^1(D)}\quad\text{for any }u\in H^1(D).
\end{equation} 
Here $p$ depends on $N$ only and $C_1$ depends on $p$, $N$ and the cone $\mathcal{C}$ only.

\end{lem}

\proof{.} See for instance the book by Adams, \cite[Theorem~5.4]{Ada}.\cvd 

\begin{oss}
We notice that if $D$ is a bounded open set such that, for some constants $p>2$ and $C_1>0$, \eqref{SobolevD} holds,
then the immersion of $H^1(D)$ into $L^2(D)$ is compact.
\end{oss}

\begin{lem}\label{secondSobolevlemma}
Let $D$ be a bounded open set satisfying the following condition. There exist constants $p>2$ and $\tilde{C}>0$ and an integer $M$ such that $D$ may be covered, up to a set of measure zero, by the union of $M$ open subsets $D_1,\ldots,D_M$ such that
$$\|u\|_{L^p(D_i)}\leq \tilde{C}\|u\|_{H^1(D_i)}\quad\text{for any }u\in H^1(D_i)\text{ and any }i=1,\ldots,M.$$ 

Then there exists a constant $C_1>0$ such that
$$\|u\|_{L^p(D)}\leq C_1\|u\|_{H^1(D)}\quad\text{for any }u\in H^1(D).$$ 
Here $C_1$ depends on $\tilde{C}$ and $M$ only.

\end{lem}

\proof{.} Let $u\in H^1(D)$. We have that
$$\|u\|^p_{L^p(D)}\leq \sum_{i=1}^M\|u_i\|^p_{L^p(D_i)}\leq \sum_{i=1}^M
\tilde{C}^p\|u\|^p_{H^1(D_i)}\leq M\tilde{C}^p \|u\|^p_{H^1(D)},$$
therefore it is enough to choose $C_1=M^{1/p}\tilde{C}$.\cvd


\subsection{The scattering case}\label{scatsubsec}

Let us consider the following scattering problem.
Let $\{K_n\}_{n\in\mathbb{N}}$ be a sequence of compact sets contained in
$\overline{B_R}$ for some $R>0$. 
Let, for any $n\in\mathbb{N}$,
$A_n=H^1(B_{R+1}\backslash K_n)$ and $\tilde{A}_n=H^1(B_{R+1}\backslash \partial K_n)$.

Let us assume that $\mathbb{R}^N\backslash K_n$ is connected and that
 the immersion of $H^1(B_{R+1}\backslash K_n)$ into $L^2(B_{R+1}\backslash K_n)$ is compact for any $n\in\mathbb{N}$. Let $u^i$ be an entire solution of the Helmholtz equation
in $\mathbb{R}^N$. Let us remark that here and in what follows we may replace, with the obvious modifications, $u^i$ with a solution of the Helmholtz equation in $\mathbb{R}^N\backslash\{x_0\}$, with $\|x_0\|>R+1$, and we may also replace $R+1$ with any $R_1>R$. Then for any $n\in\mathbb{N}$ there exists a unique weak solution to the following scattering problem
\begin{equation}\label{unscateq}
\left\{\begin{array}{ll}
\Delta u_n+k^2u_n=0 &\text{in }\mathbb{R}^N\backslash K_n\\
u_n=u^i+u^s_n &\text{in }\mathbb{R}^N\backslash K_n\\
\nabla u_n\cdot \nu=0 &\text{on }\partial K_n\\
\lim_{r\to+\infty}r^{(N-1)/2}\left(\frac{\partial u^s_n}{\partial r}-\rmi ku^s_n\right)=0 &r=\|x\|.
\end{array}\right.
\end{equation}
For existence and uniqueness see for instance \cite{Wil}.

We assume that $K_n$ converges, as $n\to\infty$, to $K\subset\overline{B_R}$ in the Hausdorff distance, where $K$ is compact and such that $\mathbb{R}^N\backslash K$ is connected.

We also assume that, as $n\to\infty$, either $A_n$ converges to $A=H^1(B_{R+1}\backslash K)$ in the sense of Mosco,
or $\tilde{A}_n$ converges to $\tilde{A}=H^1(B_{R+1}\backslash \tilde{K})$ in the sense of Mosco
and $\partial K_n$ converges to $\tilde{K}$ in the Hausdorff distance, for some $\tilde{K}\subset\overline{B_R}$, $\tilde{K}$ compact.


Let us further assume that there exist constants $p>2$ and $C_1>0$ such that for any $n\in\mathbb{N}$ we have
$$\|u\|_{L^p(B_{R+1}\backslash K_n)}\leq C_1\|u\|_{H^1(B_{R+1}\backslash K_n)}\quad\text{for any }u\in H^1(B_{R+1}\backslash K_n).$$ 

We begin with the following lemma.

\begin{lem}\label{scatlem}
Let us assume that, for some positive constant $C$, we have
\begin{equation}\label{boundedness3}
\|u_n\|_{L^2(B_{R+1})}\leq C\quad\text{for any }n\in\mathbb{N}.
\end{equation}

Then $u_n$ converges to a function $u$ strongly in $L^2(B_r)$ for any $r>0$, with $u$ solving
\begin{equation}\label{uscateq}
\left\{\begin{array}{ll}
\Delta u+k^2u=0 &\text{in }\mathbb{R}^N\backslash K\\
u=u^i+u^s &\text{in }\mathbb{R}^N\backslash K\\
\nabla u\cdot \nu=0 &\text{on }\partial K\\
\lim_{r\to+\infty}r^{(N-1)/2}\left(\frac{\partial u^s}{\partial r}-\rmi ku^s\right)=0 &r=\|x\|.
\end{array}\right.
\end{equation}
\end{lem}

\proof{.} By Lemma~3.1 in \cite{Ron03}, we have that, up to a subsequence,
$u_n$ converges to a function $u$ uniformly on compact subsets of $\mathbb{R}^N\backslash K$,
with $u$ solving
$$
\left\{\begin{array}{ll}
\Delta u+k^2u=0 &\text{in }\mathbb{R}^N\backslash K\\
u=u^i+u^s &\text{in }\mathbb{R}^N\backslash K\\
\lim_{r\to+\infty}r^{(N-1)/2}\left(\frac{\partial u^s}{\partial r}-\rmi ku^s\right)=0 &r=\|x\|.
\end{array}\right.
$$
By Proposition~\ref{firstprop} or Proposition~\ref{firstpropmodified}, without loss of generality we may assume that $u_n$ converges
to $u$ weakly in $L^2(B_{R+1})$ and $u$ also satisfies the boundary condition
$$\nabla u\cdot \nu=0\quad\text{on }\partial K,$$
that is $u$ solves \eqref{uscateq}. Since \eqref{uscateq} has at most one solution, we have that the whole sequence $u_n$ converges to $u$ uniformly on compact subsets of $\mathbb{R}^N\backslash K$.

It remains to prove that $u_n$ converges to $u$ in $L^2(B_r)$ for a fixed $r$, $R<r<R+1$.
By the reasonings used in the proof of Proposition~\ref{firstprop}, we have that there exists a constant $\tilde{C}>0$ such that
\begin{equation}
\|(u_n,\nabla u_n)\|_{L^2(B_r,\mathbb{R}^{N+1})}\leq \tilde{C}\quad\text{for any }n\in \mathbb{N}.
\end{equation}
Then we can conclude the proof by using Proposition~\ref{secondprop}.\cvd

\bigskip

In the next proposition we wish to drop the assumption \eqref{boundedness3}.

\begin{prop}\label{convergenceprop}
Let $\{K_n\}_{n\in\mathbb{N}}$ be a sequence of compact sets contained in
$\overline{B_R}$ for some $R>0$. Let us assume that $\mathbb{R}^N\backslash K_n$ is connected.

We assume that $K_n$ converges, as $n\to\infty$, to $K\subset\overline{B_R}$ in the Hausdorff distance, where $K$ is compact and such that $\mathbb{R}^N\backslash K$ is connected.

We also assume that, as $n\to\infty$, either $A_n$ converges to $A=H^1(B_{R+1}\backslash K)$ in the sense of Mosco,
or $\tilde{A}_n$ converges to $\tilde{A}=H^1(B_{R+1}\backslash \tilde{K})$ in the sense of Mosco
and $\partial K_n$ converges to $\tilde{K}$ in the Hausdorff distance, for some $\tilde{K}\subset\overline{B_R}$, $\tilde{K}$ compact.


Let us further assume that there exist constants $p>2$ and $C_1>0$ such that for any $n\in\mathbb{N}$ we have
$$\|u\|_{L^p(B_{R+1}\backslash K_n)}\leq C_1\|u\|_{H^1(B_{R+1}\backslash K_n)}\quad\text{for any }u\in H^1(B_{R+1}\backslash K_n).$$ 

Let $u_n$ solve \eqref{unscateq}, for any $n\in\mathbb{N}$.
Then $u_n$ converges to a function $u$ strongly in $L^2(B_r)$ for any $r>0$, with $u$ solving
\eqref{uscateq}.
\end{prop}

\proof{.} Let $a_n=\|u_n\|_{L^2(B_{R+1})}$. If $\{a_n\}_{n\in\mathbb{N}}$ is bounded, then the conclusion follows from the previous lemma.

By contradiction, let us assume that $\lim_{n\to\infty} a_n=+\infty$, possibly by passing to a subsequence.
Let us consider $v_n=u_n/a_n$. We have that
$$\|v_n\|_{L^2(B_{R+1})}=1.$$
Therefore $v_n$, up to a subsequence, converges to a function $v$ strongly in $L^2$ on any compact subset of $\mathbb{R}^N$.
The function $v$ satisfies
\begin{equation}
\left\{\begin{array}{ll}
\Delta v+k^2v=0 &\text{in }\mathbb{R}^N\backslash K\\
\nabla v\cdot \nu=0 &\text{on }\partial K.
\end{array}\right.
\end{equation}
Clearly we also have that $\|v\|_{L^2(B_{R+1})}=1$.

On the other hand, we have that $\|u^s_n/a_n\|_{L^2(B_{R+1})}$, $n\in\mathbb{N}$, is uniformly bounded. Therefore, again up to a subsequence, $u^s_n/a_n$ converges, as $n\to\infty$, to a function $w$ strongly in $L^2$ on any compact subset of $\mathbb{R}^N\backslash K$.
Such a function $w$ satisfies
\begin{equation}
\left\{\begin{array}{ll}
\Delta w+k^2w=0 &\text{in }\mathbb{R}^N\backslash K\\
\lim_{r\to+\infty}r^{(N-1)/2}\left(\frac{\partial w}{\partial r}-\rmi kw\right)=0 &r=\|x\|.
\end{array}\right.
\end{equation}
Since $v_n=u^i/a_n+u^s_n/a_n$, we may immediately conclude that, outside $K$, we have $v=w$.
That is $v$ solves
\begin{equation}
\left\{\begin{array}{ll}
\Delta v+k^2v=0 &\text{in }\mathbb{R}^N\backslash K\\
\nabla v\cdot \nu=0 &\text{on }\partial K\\
\lim_{r\to+\infty}r^{(N-1)/2}\left(\frac{\partial v}{\partial r}-\rmi kv\right)=0 &r=\|x\|.
\end{array}\right.
\end{equation}
By uniqueness, and the fact that $\mathbb{R}^N\backslash K$ is connected, we may conclude that
$v$ is identically zero, which leads to a contradiction since
$\|v\|_{L^2(B_{R+1})}$ should be equal to $1$.\cvd

\begin{oss}\label{variableremark}
We observe that all the results in this subsection remain valid if we let the incident field $u^i$ depend on $n$. In particular we may have that the wavenumber $k$ and the direction of propagation $d$ depend on $n$, with the assumption that $k_n>0$ converges to a real number
$k_{\infty}> 0$ and $d_n$ converges to a direction $d_{\infty}$, therefore $u^i_n$ converges to
the incident field $u^i_{\infty}(x)=\rme^{\rmi k_{\infty}x\cdot d_{\infty}}$, $x\in\mathbb{R}^N$.

We further notice that, if $N\geq 3$, then we may allow $k_{\infty}$ to be even equal to zero. In this case we need to replace the 
Sommerfeld radiation condition with the following condition at infinity for the scattered field
$$
u_s=o(1)\quad\text{as }r=\|x\|\to+\infty.
$$
We refer to \cite[Lemma~3.1]{Ron03} for further details. For $N=2$, instead, we need to require $k_{\infty}>0$, see for instance the discussion in \cite{Kre} for the low wavenumber asymptotics of scattering problems in dimension $2$.
\end{oss}

\section{Uniform decay property for scattered fields}\label{scatsec}

We begin by defining suitable classes of admissible scatterers.

\begin{defin}
We say that $\tilde{\mathcal{A}}$ is an \emph{admissible class of scatterers} in $\mathbb{R}^N$ if the following properties holds.

\begin{enumerate}[i)]
\item\label{uniformboundedness}
There exists a constant $R>0$ such that any $K\in \tilde{\mathcal{A}}$ is a compact set contained in $\overline{B_R}\subset\mathbb{R}^N$.
\item\label{connectedness} For any $K\in\tilde{\mathcal{A}}$, we have that $\mathbb{R}^N\backslash K$ is connected. Furthermore, if $K$ is the limit in the Hausdorff distance of a sequence of elements of the class
$\tilde{\mathcal{A}}$, then $\mathbb{R}^N\backslash K$ is connected.
\item\label{Moscocompactness} $\tilde{\mathcal{A}}$ satisfies one of the two following conditions.
Either for any  sequence $K_n\in\tilde{\mathcal{A}}$, $n\in\mathbb{N}$,
converging in the Hausdorff distance to $K$, we have that $H^1(B_{R+1}\backslash K_n)$ converges to $H^1(B_{R+1}\backslash K)$ in the sense of Mosco.
Or for any  sequence $K_n\in\tilde{\mathcal{A}}$, $n\in\mathbb{N}$, such that $\partial K_n$
converges in the Hausdorff distance to $\tilde{K}$, we have that $H^1(B_{R+1}\backslash \partial K)$ converges to $H^1(B_{R+1}\backslash \tilde{K})$ in the sense of Mosco.
\item\label{Sobolevcondition} There exist constants $p>2$ and $C_1>0$ such that for any $K\in\tilde{\mathcal{A}}$ we have
$$\|u\|_{L^p(B_{R+1}\backslash K)}\leq C_1\|u\|_{H^1(B_{R+1}\backslash K)}\quad\text{for any }u\in H^1(B_{R+1}\backslash K).$$
\end{enumerate}
\end{defin}

Let us then fix an admissible class of scatterers $\tilde{\mathcal{A}}$.
By condition~\ref{connectedness}), $\mathbb{R}^N\backslash K$ is connected, for any $K\in\tilde{\mathcal{A}}$. 
We also recall that condition~\ref{Sobolevcondition}) above implies that
the immersion of $H^1(B_{R+1}\backslash K)$ into $L^2(B_{R+1}\backslash K)$ is compact, for any $K\in\tilde{\mathcal{A}}$.

Let us fix constants $0<\underline{k}<\overline{k}$ and let us denote, for any $N\geq 2$,
\begin{equation}\label{I_Ndefin}
I_N=\left\{\begin{array}{ll}
[\underline{k},\overline{k}]&\text{if }N=2,\\
(0,\overline{k}]&\text{if }N\geq 3.
\end{array}\right.
\end{equation}

Fixed $K\in\tilde{\mathcal{A}}$, for a fixed wavenumber $k\in I_N$ and a fixed direction of propagation $d\in\mathbb{S}^{N-1}$, let the incident field $u^i$ be the corresponding plane wave, that is $u^i(x)=\rme^{\rmi kx\cdot d}$, $x\in\mathbb{R}^N$. Then, let $u_{K,k,d}$ be the solution to \eqref{uscateq} and $u^s_{K,k,d}$ be its corresponding scattered field.

We begin by stating the following uniform boundedness of solutions.

\begin{prop}\label{boundednesssolut}
Let us fix constants $0<\underline{k}<\overline{k}$ and let $I_N$ be defined as in \eqref{I_Ndefin}.
Let $\tilde{\mathcal{A}}$ be an admissible class of scatterers in $\mathbb{R}^N$.

Fixed $K\in\tilde{\mathcal{A}}$, $k\in I_N$, and $d\in\mathbb{S}^{N-1}$, let 
$u^i(x)=\rme^{\rmi kx\cdot d}$, $x\in\mathbb{R}^N$, $u_{K,k,d}$ be the solution to \eqref{uscateq} and $u^s_{K,k,d}$ be its corresponding scattered field.

Then there exists a constant $E$, depending on $\tilde{\mathcal{A}}$ and
$I_N$ only, such that
\begin{equation}\label{uniformbound}
\|u_{K,k,d}\|_{L^2(B_{R+1})}\leq E\quad\text{for any }K\in \tilde{\mathcal{A}},\text{ any }k\in I_N,\text{ and any }d\in\mathbb{S}^{N-1}.
\end{equation}
\end{prop}

\proof{.} We sketch the proof. We argue by contradiction. Let us assume that there exist, for any $n\in\mathbb{N}$,
$K_n\in\tilde{\mathcal{A}}$,
$k_n\in I_N$, and $d_n\in\mathbb{S}^{N-1}$ such that
$$\|u_{K_n,k_n,d_n}\|_{L^2(B_{R+1})}\geq n.$$
Then we define $v_n=u_{K_n,k_n,d_n}/\|u_{K_n,k_n,d_n}\|_{L^2(B_{R+1})}$, $n\in\mathbb{N}$. 

Then we obtain a contradiction arguing in an analogous manner as in Proposition~\ref{convergenceprop} and using Remark~\ref{variableremark}.\cvd

\bigskip

As a corollary, as in Lemma~3.2 in \cite{Ron03}, we may prove the following uniform decay property.

\begin{cor}\label{uniformdecay}
Under the same assumptions and notation of Proposition~\textnormal{\ref{boundednesssolut}},
there exists a constant $E_1$, depending on the constant  $E$ in \eqref{uniformbound}, $I_N$, $R$ and $N$ only,
such that for any $K\in\tilde{\mathcal{A}}$, any $k\in I_N$, and any $d\in\mathbb{S}^{N-1}$
we have
$$|u^s_{K,k,d}(x)|\leq E_1\|x\|^{-(N-1)/2}\quad\text{for any }x\in\mathbb{R}^N\text{ so that }\|x\|\geq R+2.$$
\end{cor}


In the remaining part of this section
we wish to discuss Assumptions~\ref{uniformboundedness})-\ref{connectedness})-\ref{Moscocompactness})-\ref{Sobolevcondition}) defined above.

Assumption~\ref{uniformboundedness}) is self-explanatory. We notice however that is guarantees 
that $\tilde{\mathcal{A}}$ is relatively compact with respect to the Hausdorff distance.

A sufficient condition for Assumption~\ref{connectedness})
to hold is the following, see \cite[Lemma~2.5]{Ron03}.

\medskip

\noindent
\textbf{Assumption~A (uniform exterior connectedness)}
Let $\delta:(0,+\infty)\to (0,+\infty)$ be a nondecreasing left-continuous function.
We assume that for any $K\in\tilde{\mathcal{A}}$,
for any $t>0$,
for any two points $x_1$, $x_2\in\mathbb{R}^N$ so that
$B_t(x_1)$ and $B_t(x_2)$ are contained in $\mathbb{R^N}\backslash K$, and for any $s$, $0<s<\delta(t)$,
then we can find a smooth (for instance $C^1$) curve $\gamma$ connecting
$x_1$ to $x_2$ so that $B_{s}(\gamma)$ is contained in
$\mathbb{R^N}\backslash K$ as well.

\bigskip

Let us notice that such an Assumption~A is closed under convergence in the Hausdorff distance.

Assumption~\ref{Sobolevcondition}) has been already discussed in Lemma~\ref{conelemma}
and Lemma~\ref{secondSobolevlemma}. Let us further notice that such an assumption is closed under convergence in the sense of Mosco of the corresponding Sobolev spaces. In particular, it holds for any $K$ which is  the limit in the Hausdorff distance of a sequence of elements of the class
$\tilde{\mathcal{A}}$, provided Assumptions~\ref{uniformboundedness}) and  \ref{Moscocompactness}) hold true, possibly by using Proposition~\ref{convboundary}.

About Assumption~\ref{Moscocompactness}) we state three sufficient conditions.
The first one is proved \cite{Cha-Dov} and it holds only for $N=2$, the second is taken from \cite{Gia} and it is valid in any dimension $N\geq 2$. The third one is new, and it will be proven here, and it uses a class which is a generalization of one developed in \cite{Ron06}.
Let us finally remark that in what follows Assumption~\ref{uniformboundedness}) will be always tacitly assumed.

\medskip

\noindent
\textbf{Assumption~B \cite{Cha-Dov}}
Let us assume that $N=2$. Let us assume that there exist constants $M\in\mathbb{N}$ and $C>0$ such that for any $K\in\tilde{\mathcal{A}}$ we have that the number of connected components of $\partial K$ is bounded by $M$ and $\mathcal{H}^1(\partial K)\leq C$.

\bigskip

We notice that, under Assumption~B, for any $K \in \tilde{\mathcal{A}}$ we have that the number of connected components of $K$ is bounded by $M$. Moreover, if $K$ is the limit in the Hausdorff distance of a sequence $K_n$, $n\in\mathbb{N}$, of sets belonging to $\tilde{\mathcal{A}}$, we have that the number of connected components of $K$ is bounded by $M$. Furthermore, without loss of generality, we may assume that $\partial K_n$ converges, in the Hausdorff distance, to a compact set $\tilde{K}$. Since $\partial K\subset \tilde{K}\subset K$, we deduce by a general version of Gol\c ab's Theorem, see for instance Corollary~3.3 in \cite{DM-Toa}, that
$H^1(\partial K)\leq H^1(\tilde{K})\leq C$ as well. We finally point out that the number of connected components of $\tilde{K}$ is bounded by $M$, whereas the same may not be true for $\partial K$.
However, if Assumption~\ref{connectedness}) holds, then $\mathbb{R}^2\backslash K$
is connected, therefore 
$K$ and $\partial K$ have the 
same number of connected components, so also the number of connected components of $\partial K$ is bounded by $M$. In other words, if Assumption~\ref{connectedness}) holds, then
Assumption~B is closed under convergence in the Hausdorff distance.

\medskip

\noindent
\textbf{Assumption~C \cite{Gia}}
There exist a fixed cone $\mathcal{C}$ in $\mathbb{R}^{N-1}$ and positive constants $\delta$, $L_1$ and $L_2$ such that any $K\in\tilde{\mathcal{A}}$ satisfies the following condition.

For any $x\in\partial K$ there exists a bi-Lipschitz function $\Phi_x: B_{\delta}(x)\to\mathbb{R}^N$ such that
\begin{enumerate}[C1)]
\item for any $z_1$, $z_2\in B_{\delta}(x)$ we have
$$L_1\|z_1-z_2\|\leq\|\Phi_x(z_1)-\Phi_x(z_2)\|\leq L_2\|z_1-z_2\|;$$
\item $\Phi_x(x)=0$ and
$\Phi_x(\partial K\cap B_{\delta}(x))\subset \pi=\{y\in\mathbb{R}^N:\ y_N=0\}$;
\item for any $y\in\partial K\cap B_{\delta/2}(x)$ we have
$$\Phi_x(y)\in \overline{\mathcal{C}_y}\subset \Phi_x(\partial K\cap B_{\delta}(x)),$$
where $\mathcal{C}_y$ is obtained by the cone $\mathcal{C}$ through a rigid change of coordinates.
\end{enumerate}

\bigskip

It is not difficult to see that Assumption~C is closed under convergence in the Hausdorff distance.
With the following lemma we show that Assumption~C guarantees that not only Assumption~\ref{Moscocompactness}) is satisfied but also Assumption~\ref{Sobolevcondition}) is.

\begin{lem}
Let us assume that $\tilde{\mathcal{A}}$ satisfies Assumption~\textnormal{C}. Then there exist constants $p>2$ and $C_1>0$ such that for any $K\in\tilde{\mathcal{A}}$ we have
$$\|u\|_{L^p(B_{R+1}\backslash \partial K)}\leq C_1\|u\|_{H^1(B_{R+1}\backslash\partial K)}\quad\text{for any }u\in H^1(B_{R+1}\backslash\partial K),$$
hence the same property is satisfied by $B_{R+1}\backslash K$.
\end{lem}

\proof{.}
We sketch the proof. Without loss of generality, let $\delta$ be such that $0<\delta\leq 1/2$.
We can find positive constants $\delta_1$, $\delta_2$ such that for any $x\in\partial  K$ we have that $B_{\delta_2}\subset \Phi_x(B_{\delta}(x))$ and $B_{\delta_1}(x)\subset
\Phi_x^{-1}(B_{\delta_2})$. Clearly $\delta_1$ and $\delta_2$ depend on $\delta$, $L_1$ and $L_2$ only.

Then we obtain that $\Phi_x^{-1}(B_{\delta_2})\backslash\partial K$ is contained, up to a set of measure zero, in the set
 $U_x=\Phi_x^{-1}(B_{\delta_2}\backslash \pi)$. Since $U_x$ is the image through a bi-Lipschitz 
map of $B_{\delta_2}\backslash \pi$, we have that $U_x$ satisfies, for some $p>2$ and $C>0$ depending on $\delta$, $L_1$ and $L_2$ only, the following Sobolev inequality
$$\|u\|_{L^p(U_x)}\leq C\|u\|_{H^1(U_x)}\quad\text{for any }u\in H^1(U_x).$$

We have that $\overline{B_{\delta_1/4}(\partial K)}$ is contained in $\bigcup_{x\in\overline{B_{\delta_1/4}(\partial K)}}B_{\delta_1/4}(x)$.
We can find a finite number of points $z_i\in \overline{B_{\delta_1/4}(\partial K)}$, $i=1,\ldots,m_1$, such that $\overline{B_{\delta_1/4}(\partial K)}\subset \bigcup_{i=1}^{m_1}B_{\delta_1/4}(z_i)$. With a simple construction, it is possible to choose $m_1$ depending on $\delta_1$ and $R$ only, for instance by taking points such that $B_{\delta_1/8}(z_i)\cap B_{\delta_1/8}(z_j)$ is empty for $i\neq j$. Then, for any $i=1,\ldots,m_1$ we can find $x_i\in \partial K$ such that
$B_{\delta_1/4}(z_i)\subset B_{\delta_1}(x_i)$, therefore $\overline{B_{\delta_1/4}(\partial K)}\subset \bigcup_{i=1}^{m_1}B_{\delta_1}(x_i)$.

Then consider the set $A=\overline{B_{R+1/2}}\backslash B_{\delta_1/4}(\partial K)$. Again, we can find points $y_j\in A$,  $j=1,\ldots,m_2$, such that $A\subset \bigcup_{j=1}^{m_2}B_{\delta_1/4}(y_j)$.
Again, with the same kind of construction, it is possible to choose $m_2$ depending on $\delta_1$ and $R$ only. We notice that $B_{\delta_1/4}(y_j)\subset B_{R+1}\backslash \partial K$ for any $j=1,\ldots,m_2$.

Therefore, $B_{R+1}\backslash K$ is contained, up to a set of measure zero, in the following union
$$\left(B_{R+1}\backslash \overline{B_{R+1/2}}\right)\cup\left(\bigcup_{i=1}^{m_1}U_{x_i}\right)\cup\left(\bigcup_{j=1}^{m_2}B_{\delta_1/4}(y_j)\right).$$
The conclusion immediately follows by Lemma~\ref{secondSobolevlemma}.\cvd

\bigskip

The third sufficient condition is a generalization of arguments developed in \cite{Ron06}. We need to fix some preliminary notation and prove some lemmas.

Let us fix a bounded open set $\Omega$.
Let $K\subset\overline{\Omega}$ be a compact subset of $\mathbb{R}^N$. We say that $K$ is a \emph{Lipschitz hypersurface}, with or without boundary, with positive constants $r$ and $L$ if the following holds.

For any $x\in K$ there exists a bi-Lipschitz function $\Phi_x: B_{r}(x)\to\mathbb{R}^N$ such that
\begin{enumerate}[a)]
\item\label{conditiona} for any $z_1$, $z_2\in B_{r}(x)$ we have
$$L^{-1}\|z_1-z_2\|\leq\|\Phi_x(z_1)-\Phi_x(z_2)\|\leq L\|z_1-z_2\|;$$
\item\label{conditionb} $\Phi_x(x)=0$ and
$\Phi_x(K\cap B_{r}(x))\subset \pi=\{y\in\mathbb{R}^N:\ y_N=0\}$;
\newcounter{enumi_saved}
\setcounter{enumi_saved}{\value{enumi}}
\end{enumerate}

We say that $x\in K$ belongs to the interior of $K$ if there exists $\delta$, $0<\delta\leq r$, such that $B_{\delta}(0)\cap \pi \subset \Phi_x(K\cap B_{r}(x))$.
Otherwise we say that $x$ belongs to the boundary of $K$. We remark that the boundary of $K$ might be empty. Further we assume that

\begin{enumerate}[a)]\setcounter{enumi}{\value{enumi_saved}}
\item\label{conditionc} for any $x$ belonging to the boundary of $K$, we have that
$$\Phi_x(K\cap B_{r}(x))=\Phi_x(B_{r}(x))\cap \pi^+$$
where $\pi^+=\{y\in\mathbb{R}^N:\ y_N=0,\ y_{N-1}\geq 0\}$.
\end{enumerate}

Let us notice that, by compactness, such an assumption is enough to guarantee that $\mathcal{H}^{N-1}(K)$ is bounded, hence $|K|=0$. In particular, if $\Omega\subset \overline{B_R}$ for some $R>0$, then $\mathcal{H}^{N-1}(K)$ is bounded by a constant depending on $R$, $r$ and $L$ only. Furthermore, the boundary of $K$ has $\mathcal{H}^{N-2}$ measure bounded by a constant again depending on $R$, $r$ and $L$ only.

Moreover, $K$ has a finite number of connected components, again bounded a constant depending on $R$, $r$ and $L$ only, and the distance between two different connected components of $K$ is bounded from below by a positive constant depending on $r$ and $L$ only.

We begin with the following lemma.

\begin{lem}\label{singlehyperlemma}
Let us fix a bounded open set $\Omega$ and positive constants $r$ and $L$.
Let $\mathcal{B}=\mathcal{B}(r,L)$ be the class of compact sets $K\subset\overline{\Omega}$ such that
$K$ is a Lipschitz hypersurface with constants $r$ and $L$.

Then $\mathcal{B}$ is closed under the convergence in the Hausdorff distance, that is if $\{K_n\}_{n\in\mathbb{N}}$ is a sequence in $\mathcal{B}$, then up to a subsequence $K_n$ converges, as $n\to\infty$, to a set $K\in\mathcal{B}$ in the Hausdorff distance. Furthermore, the boundary of $K_n$ converges to the boundary of $K$ in the Hausdorff distance.
\end{lem}

\proof{.} Without loss of generality, up to a subsequence, we may assume that $K_n$ and their boundaries converge to a compact set $K$ and a compact set $H$, respectively.

Let $x\in K$. Then there exists a sequence $x_n\in K_n$ such that $\lim_{n\to\infty}x_n=x$. Again up to subsequences, we have that $\Phi_{x_n}^n$ converges to a function $\Phi_x:B_{r}(x)\to\mathbb{R}^N$ satisfying Condition~\ref{conditiona}) above. Clearly, since $\Phi_{x_n}^n(x_n)=0$, we have that $\Phi_{x}(x)=0$ as well. It is not difficult to show that Condition~\ref{conditionb}) is also satisfied.

As far as Condition~\ref{conditionc}) is concerned, we first prove that $H$ coincide with the boundary of $K$.
Let $x\in H$, then there exists $x_n$ belonging to the boundary of $K_n$ such that $\lim_{n\to\infty}x_n=x$.
Again up to subsequences, we have that $\Phi_{x_n}^n$ converges to a function $\Phi_x:B_{r}(x)\to\mathbb{R}^N$ satisfying Condition~\ref{conditiona}) above. Again it is not difficult to prove that $\Phi_x$ satisfies both Condition~\ref{conditionb}) and Condition~\ref{conditionc}). Therefore, $H$ is contained in the boundary of $K$. We need to prove that the boundary of $K$ is contained in $H$. By contradiction, let us assume that there exists $x$ belonging to 
the boundary of $K$ such that $x$ does not belong to $H$. For some positive constant $c$ and for any $n$ large enough, we have that the distance of $x$ from the boundary of $K_n$ is greater than $c$. If $x_n\in K_n$ converges to $x$, we have that $x_n$ has a distance from the boundary of $K_n$ greater than $c/2$, for any $n$ large enough.
Therefore there exists a positive constant $c_1$ such that for any $n$ large enough we have
$B_{c_1}(0)\cap \pi \subset \Phi^n_{x_n}(K_n\cap B_{r}(x_n))$. Passing to the limit, we obtain that
$B_{c_1}(0)\cap \pi \subset \Phi_{x}(K\cap B_{r}(x))$ as well, therefore $x$ belongs to the interior of $K$. Such a contradiction concludes the proof.\cvd

\begin{defin}\label{classdefin}
Let us fix a bounded open set $\Omega$ and positive constants $r$ and $L$. We assume that $\partial\Omega$ consists of a finite number of hypersurfaces without boundary belonging to $\mathcal{B}(r,L)$

Let us also fix $\omega:(0,+\infty)\to (0,+\infty)$ a nondecreasing left-continuous function.

We say that a compact set $K\subset \overline{\Omega}$ belongs to the class $\tilde{\mathcal{B}}=
\tilde{\mathcal{B}}(r,L,\omega)$ if $\hat{K}=K\cup\partial\Omega$ satisfies the following conditions.

\begin{enumerate}[1)]
\item\label{1} $\hat{K}=\bigcup_{i=1}^M K^i$ where $K^i\in\mathcal{B}(r,L)$ for any $i=1,\ldots,M$.
\item\label{3} For any $i\in \{1,\ldots,M\}$, and any $x\in K^i$, if its distance from the boundary of $K^i$ is $\delta>0$, then the distance of $x$ from the union of $K^j$, with $j\neq i$, is greater than or equal to $\omega(\delta)$.
\end{enumerate}
\end{defin}

Let us notice that in the previous definition the number $M$ may depend on $K$. However, there exists an integer $M_0$, depending on the diameter of $\Omega$, $r$, $L$ and $\omega$ only, such that $M\leq M_0$ for any $K\in \tilde{\mathcal{B}}$.
As before, we obtain that $\mathcal{H}^{N-1}(\hat{K})$ is bounded, hence $|\hat{K}|=0$. In particular, if $\Omega\subset \overline{B_R}$ for some $R>0$, then $\mathcal{H}^{N-1}(\hat{K})$ is bounded by a constant depending on $R$, $r$, $L$ and $M_0$ only. Furthermore, if we set as the boundary of $\hat{K}$ the union of the boundaries of $K^i$, $i=1,\ldots,M$, then
the boundary of $\hat{K}$ has $\mathcal{H}^{N-2}$ measure bounded by a constant again depending on $R$, $r$, $L$ and $M_0$ only.

We also remark that, by Condition~\ref{3}), we have that $K^i\cap K^j$ is contained in the intersection of the boundaries of $K^i$ and $K^j$, for any $i\neq j$.

We prove the analogous of Lemma~\ref{singlehyperlemma} for the class $\tilde{\mathcal{B}}$.

\begin{lem}\label{multiplehyperlemma}
Under the previous notation and assumptions, we have that $\tilde{\mathcal{B}}=\tilde{\mathcal{B}}(r,L,\omega)$ as in Definition~\textnormal{\ref{classdefin}} is closed under the convergence in the Hausdorff distance, that is if $\{K_n\}_{n\in\mathbb{N}}$ is a sequence in $\tilde{\mathcal{B}}$, then up to a subsequence $K_n$ converges, as $n\to\infty$, to a set $K\in\tilde{\mathcal{B}}$ in the Hausdorff distance. Clearly also $\hat{K_n}$ converges to $\hat{K}$ in the Hausdorff distance. Furthermore, if $\tilde{K}_n$, $n\in\mathbb{N}$, is the corresponding sequence of boundaries of $\hat{K}_n$ and $\tilde{K}$ is the boundary of $\hat{K}$, then
$\tilde{K}_n$ converges to $\tilde{K}$  in the Hausdorff distance.
\end{lem}

\proof{.} Up to a subsequence, we may assume that as $n\to\infty$ $K_n$ converges to a set $K$, hence $\hat{K_n}$ converges to $\hat{K}=K\cup\partial\Omega$, and $\tilde{K}_n$ converges to a set $H$. Moreover, we may assume that $M(n)=M$ for any $n\in\mathbb{N}$ and,
by Lemma~\ref{singlehyperlemma}, that for any $i=1,\ldots,M$ $K_n^i$ converges to a set $K^i\in \mathcal{B}=\mathcal{B}(r,L)$
and the boundary of $K^i_n$ converges to the boundary of $K^i$. It is not difficult then to show that $\hat{K}=\bigcup_{i=1}^M K^i$ and that $H$ is the union of the boundaries of $K^i$, $i=1,\ldots,M$. Therefore $K$ satisfies Condition~\ref{1}) above.

We now deal with Condition~\ref{3}).
Let us take $x$ belonging to the interior of $K^i$ and let $\delta>0$ be its distance from the boundary of $K^i$. Let $x_n\in K^i_n$ converge to $x$. 
For any $\varepsilon$, $0<\varepsilon<\delta/2$, there exists $\overline{n}=\overline{n}(\varepsilon)$ such that for any $n\geq \overline{n}$ we have that the distance of $x_n$ from the boundary of $K^i_n$ is greater than $\delta-\varepsilon$.
We obtain that the distance of $x_n$ from the union of $K^j_n$, with $j\neq i$, is greater than or equal to $\omega(\delta-\varepsilon)$, for any $n\geq \overline{n}$. Provided $0<\varepsilon<\omega(\delta-\varepsilon)$, we can find a further $\hat{n}=\hat{n}(\varepsilon)$ such that for any $n\geq \hat{n}$ we have that the distance of $x_n$ from the union of $K^j$, with $j\neq i$, is greater than or equal to $\omega(\delta-\varepsilon)-\varepsilon$. Passing to the limit as $n\to\infty$ we obtain that 
the distance of $x$ from the union of $K^j$, with $j\neq i$, is greater than or equal to $\omega(\delta-\varepsilon)-\varepsilon$ for any $0<\varepsilon<\omega(\delta-\varepsilon)$. We then let $\varepsilon\to 0^+$ and we conclude by the 
left-continuity of $\omega$.\cvd

\bigskip

We are now in the position of stating and proving the following Mosco convergence result.

\begin{teo}\label{Moscoconvergencetheo}
Let $\Omega$ and $\tilde{\mathcal{B}}=\tilde{\mathcal{B}}(r,L,\omega)$ be as in Definition~\textnormal{\ref{classdefin}}.

Let $\{K_n\}_{n\in\mathbb{N}}$ be a sequence in $\tilde{\mathcal{B}}$ converging, as $n\to\infty$, to $K\in \tilde{\mathcal{B}}$ in the Hausdorff distance.

Then $H^1(\Omega\backslash K_n)$ converges, as $n\to\infty$, to $H^1(\Omega\backslash K)$ in the sense of Mosco.
\end{teo}

\proof{.}
Let us denote $A_n=H^1(\Omega\backslash K_n)$, $n\in\mathbb{N}$, and $A=H^1(\Omega\backslash K)$.

Since $|K|=0$, then by Lemma~\ref{primolemma}, we immediately have that $A'\subset A$. Therefore it is enough to prove that $A\subset A''$ or, in other words, that for every $\varphi\in A$ there exists $\varphi_n\in A_n$ such that $\varphi_n$ converges as $n\to\infty$ to $\varphi$ in $L^2(\Omega,\mathbb{R}^{N+1})$. We notice that it is enough to prove that for any subsequence
$A_{n_{k}}$ there exists a further subsequence $A_{n_{k_j}}$ and $\varphi_j\in A_{n_{k_j}}$
such that $\varphi_j$ converges as $j\to\infty$ to $\varphi$ in $L^2(\Omega,\mathbb{R}^{N+1})$.
Therefore during our proof we can always pass to subsequences, without loss of generality.

We may assume that $\hat{K}_n$ is converging to $\hat{K}$ and $\tilde{K}_n$ is converging to
$\tilde{K}$ in the Hausdorff distance.

Since $A''$ is closed, it is enough to prove the result for any $\varphi$ in a dense subset of $A$.
For instance, let us consider the following subset of $A$
$$\tilde{A}=\{\varphi\in H^1(\Omega\backslash K):\ \varphi\text{ is bounded and }\varphi=0\text{ in a neighbourhhod of }\tilde{K}\}.$$
We wish to show that $\tilde{A}$ is dense in $A$. By an easy truncation argument we can show that $\{\varphi\in H^1(\Omega\backslash K):\ \varphi\text{ is bounded}\}$ is dense in $A$. It is enough to show that $\tilde{A}$ is dense in this last set.
Since $\mathcal{H}^{N-2}(\tilde{K})$ is finite, then $\tilde{K}$ has zero capacitiy. Hence for any neighborhood $U$ of  $\tilde{K}$ and for any $\varepsilon>0$
there exists a function $\chi_{\varepsilon}$ such that $\chi_{\varepsilon}\in H^1(\Omega)$, $0\leq\chi_{\varepsilon}\leq 1$ almost everywhere in $\Omega$, $\chi_{\varepsilon}=1$ almost everywhere outside $U$, $\chi_{\varepsilon}=0$ almost everywhere in a neighbourhood of $\tilde{K}$, and
$$\int_{\Omega}\|\nabla \chi_{\varepsilon}\|^2\leq \varepsilon.$$
Take $\varphi\in H^1(\Omega\backslash K)$ such that $\varphi$ is bounded.
Clearly we have that $\chi_{\varepsilon}\varphi\in \tilde{A}$ and
$$\|\chi_{\varepsilon}\varphi-\varphi\|_{L^2(\Omega)}\leq \|\varphi\|_{L^2(U)},
\quad
\|\nabla (\chi_{\varepsilon}\varphi)-\nabla\varphi\|_{L^2(\Omega)}\leq  \|\nabla \varphi\|_{L^2(U)}+
\|\varphi\|_{L^{\infty}(\Omega)}\sqrt{\varepsilon}.
$$
Since $U$ and $\varepsilon$ are arbitrary, we conclude that $\tilde{A}$ is dense in $A$.

Take $\varphi\in\tilde{A}$ and let $\tilde{U}$ be an open neighborhood of $\tilde{K}$ on which $\varphi$ is zero. We can find an open subset $D$ compactly contained in $\Omega\backslash K$, a finite number of points $x_j\in \hat{K}$ and positive numbers $\delta_j$, $j=1,\ldots,m$, such that
$B_{2\delta_j}(x_j)\cap \tilde{K}=\emptyset$ for any $j=1,\ldots,m$ and
$$\overline{\Omega}\subset(\tilde{U}\cap \Omega)\cup D\cup\left(\bigcup_{j=1}^m (B_{\delta_j}(x_j)\cap\Omega)\right).$$
Moreover, for any $j=1,\ldots,m$, we have that $x_j$ belongs to the interior of $K^{i(j)}$ for some $i=i(j)\in\{1,\ldots,M\}$ and we may assume that 
$2\delta_j\leq r$ and that $B_{2\delta_j}(x_j)\cap \hat{K}=B_{2\delta_j}(x_j)\cap K^{i(j)}$.

By using a partition of unity, we may therefore reduce ourselves to the following cases. It is sufficient to consider a function $\varphi\in\tilde{A}$ that is compactly supported either in $B_{\delta_j}(x_j)$, for some $j\in\{1,\ldots,m\}$, or in $D$. In the latter case, we have that $D\subset (\Omega
\backslash K_n)$, hence $\varphi\in H^1(\Omega\backslash K_n)$, for any $n$ large enough, so the convergence is trivially proved.

It remains to prove the convergence for $\varphi\in\tilde{A}$ that is compactly supported
in $B_{\delta}(x)$ for some $0<\delta\leq r/2$ and $x\in K^i\subset\hat{K}$ such that
$B_{2\delta}(x)\cap \tilde{K}=\emptyset$ and
$B_{2\delta}(x)\cap \hat{K}=B_{2\delta}(x)\cap K^i$.

We use the reasoning developed in the proof of Theorem~4.2 in \cite{Gia}. For the convenience of the reader we repeat the construction. We suppose that $x\in\Omega$, the case when $x\in\partial\Omega$ requires just a little modification.
Possibly passing to a subsequence, let $x_n\in K_n$ converge to $x$ and
$\Phi^n_{x_n}$ converge to a function $\Phi_{x}:B_{r}(x)\to\mathbb{R}^N$.
Without loss of generality, we may assume that $B_{\delta}(x)\subset \Phi_x^{-1}(B_{r_1})$ for some positive $r_1$ such that
$B_{r_1}\subset \Phi_x(B_r(x))$. Moreover, we may also assume that $B_{r_1}\cap \pi\subset \Phi_x(B_r(x)\cap K)$ and $B_{r_1}\cap \pi\subset \Phi^n_{x_n}(B_r(x_n)\cap K_n)$ for any $n$
large enough.
Let $\psi=\varphi\circ\Phi_x^{-1}$. Then $\psi\in H^1(B_{r_1}\backslash \pi)$.
We denote by $\psi^{\pm}$ the function $\psi$ defined above or below $\pi$, that is in the halfspaces $T^{\pm}=\{y\in\mathbb{R}^N:\ \pm y_n>0\}$, respectively. Then, by an odd reflection, we may define two $H^1_0(B_{r_1})$ functions, $\tilde{\psi}^{\pm}$ such that $\tilde{\psi}^{\pm}=\psi^{\pm}$ on $T^{\pm}$. Let $\tilde{\varphi}^{\pm}=\tilde{\psi}^{\pm}\circ\Phi_x\in H^1_0(B_{\delta}(x))$.

Finally we define
$$\varphi_n=\left\{\begin{array}{ll}
\tilde{\varphi}^{+}(x) &\text{if }\Phi^n_{x_n}(x)\in T^+\\
\tilde{\varphi}^{-}(x) &\text{if }\Phi^n_{x_n}(x)\in T^-.
\end{array}\right.
$$
By construction we have that $\varphi_n\in H^1(\Omega\backslash K_n)$. Furthermore, $\varphi_n$ converges almost everywhere to $\varphi$. Since $|\varphi_n|\leq \max\{|\tilde{\varphi}^{+}|,|\tilde{\varphi}^{-}|\}$,
 by the Dominated Convergence Theorem we have that $\varphi_n$ converges to $\varphi$ in $L^2$.
The same argument holds true for the gradients, so the proof is concluded.\cvd

\bigskip

As an immediate corollary to Theorem~\ref{Moscoconvergencetheo}, we infer that the following assumption on $\tilde{\mathcal{A}}$ is a sufficient condition for Assumption~\ref{Moscocompactness}) to hold.

\medskip

\noindent
\textbf{Assumption~D}
Fixed $R>0$, let $\Omega=B_{R+1}$. Let us fix positive constants $r$ and $L$ and let $\omega:(0,+\infty)\to (0,+\infty)$ be a nondecreasing left-continuous function.

We assume that $\tilde{\mathcal{A}}$ satisfies Assumption~\ref{uniformboundedness}) and that
for any $K\in \tilde{\mathcal{A}}$ we have that $\partial K$ belongs to $\tilde{\mathcal{B}}(r,L,\omega)$.

\bigskip

We conclude the section simply by pointing out that, by Lemma~\ref{multiplehyperlemma}, we have that Assumption~D is closed 
under convergence in the Hausdorff distance.

%
%
%
%
%

\section{Approximation of sound-hard screens}\label{examplesec}

Let $K$ be a compact set 
such that $\mathbb{R}^N\backslash K$ is connected. Let us assume that
$K\subset \overline{B_R}$ for some $R>0$.

Let us denote $d:\mathbb{R}^N\to[0,+\infty)$ the function defined as follows
$$d(x)=\mathrm{dist}(x,K)\quad\text{for any }x\in\mathbb{R}^N.$$

We wish to find sufficient conditions on $K$ so that 
there exists a Lipschitz function $\tilde{d}:\mathbb{R}^N\to [0,+\infty)$ such that
the following properties are satisfied.

First, there exist constants $a$ and $b$, $0<a\leq1\leq b$, such that
$$ad(x)\leq\tilde{d}(x)\leq bd(x)\quad\text{for any }x\in\mathbb{R}^N.$$

For any $h>0$, let us call $K_h=\{x\in\mathbb{R}^N:\ \tilde{d}(x)\leq h\}$.
Second, we require that, for some constants $h_0>0$, $p>2$ and $C_1>0$, for any $h$, $0<h\leq h_0$, $\mathbb{R}^N\backslash K_h$ is connected
and 
$$\|u\|_{L^p(B_{R+1}\backslash K_h)}\leq C_1\|u\|_{H^1(B_{R+1}\backslash K_h)}\quad\text{for any }u\in H^1(B_{R+1}\backslash K_h).$$ 

Then, if we take a decreasing sequence of positive numbers $\{\varepsilon_n\}_{n\in\mathbb{N}}$ such that $\lim_n\varepsilon_n=0$ and define $K_n=K_{\varepsilon_n}$ for any $n\in\mathbb{N}$,
we immediately have that $K_n$ converges to $K$ in the Hausdorff distance and that $H^1(B_{R+1}\backslash K_n)$ converges to $H^1(B_{R+1}\backslash K)$ in the sense of Mosco. Finally, all the assumptions of Proposition~\ref{convergenceprop} are satisfied.


%

A simple sufficient condition is that $K$ is a compact convex set. In fact, clearly we have that 
$\mathbb{R}^N\backslash K$ is connected. Then,
we can take $\tilde{d}=d$ or the distance from $K$ with respect to any norm on $\mathbb{R}^N$, not only with respect to the Euclidean one. Then, for any $h>0$, $K_h=\overline{B_h(K)}$, clearly with respect to the chosen norm.
For any $h>0$, we have that $\overline{B_h(K)}$ is still a convex set, therefore $\mathbb{R}^N\backslash \overline{B_h(K)}$ is connected and also the other required properties are satisfied, for example by using Lemma~\ref{conelemma}.
Clearly we can extend this property to a set $K$ which is the union of a finite number of 
compact convex sets which are pairwise disjoint.

In what follows we investigate sufficient conditions for a set $K$ satisfying some minimal Lipschitz type regularity assumptions.
The assumptions on $K$ are the following. Let us fix positive constants $r$ and $L$.

For any $x\in\partial K$, there exists a function $\varphi:\mathbb{R}^{N-1}\to\mathbb{R}$, such that $\varphi(0)=0$ and which is Lipschitz with Lipschitz constant bounded by $L$, such that, up to a rigid change of coordinates, we have $x=0$ and
$$B_r(x)\cap \partial K\subset \{y\in B_r(x): y_N=\varphi(y')\}.$$

Let us notice that, by compactness, such an assumption is enough to guarantee that $\mathcal{H}^{N-1}(\partial K)$ is bounded, hence $|\partial K|=0$. In particular, if $K\subset \overline{B_R}$ for some $R>0$, then $\mathcal{H}^{N-1}(\partial K)$ is bounded by a constant depending on $R$, $r$ and $L$ only.

We say that $x\in \partial K$ belongs to the interior of $\partial K$ if there exists $\delta$, $0<\delta\leq r$, such that $B_{\delta}(x)\cap \partial K = \{y\in B_{\delta}(x): y_N=\varphi(y')\}$. Otherwise we say that $x$ belongs to the boundary of $\partial K$. We remark that the boundary of $\partial K$ might be empty and that, if $x\in \partial K$ belongs to the interior of $\partial K$, then $K$ may lie at most on one side of $\partial K$, that is $B_{\delta}(x)\cap K=B_{\delta}(x)\cap \partial K$, or
$B_{\delta}(x)\cap K=\{y\in B_{\delta}(x): y_N\geq \varphi(y')\}$, or 
$B_{\delta}(x)\cap K=\{y\in B_{\delta}(x): y_N\leq \varphi(y')\}$.

For any $x$ belonging to the boundary of $\partial K$, we assume  that there exists another function
$\varphi_1:\mathbb{R}^{N-2}\to\mathbb{R}$, such that $\varphi_1(0)=0$ and which is Lipschitz with Lipschitz constant bounded by $L$, such that, up to the previous rigid change of coordinates, we have $x=0$ and
$$B_r(x)\cap \partial K= \{y\in B_r(x): y_N=\varphi(y'),\ y_{N-1}\leq\varphi_1(y'')\}.$$

Let us remark that by this properties we have that $\partial K$ is a Lipschitz hypersurface as
defined in Section~\ref{scatsec}.

We notice that $K$ has a finite number of connected components, again bounded a constant depending on $R$, $r$ and $L$ only, and the distance between two different connected components of $K$ is bounded from below by a positive constant depending on $r$ and $L$ only. More precisely, we have that each connected component of $K$ is either the closure of a Lipschitz domain or a Lipschitz hypersurface with boundary. Obviously, the exterior of any connected component of $K$ is connected, since the exterior of $K$ is. Furthermore, 
the numbers of connected components of $K$ and $\partial K$coincide, that is if a connected component of $K$ is the closure of a Lipschitz domain then its boundary is connected. Here we have made use again  of the fact that the exterior of $K$ is connected. We may conclude that
the exterior of $K$ is connected if and only if the exteriors of its connected components are connected.

%

Finally, for any $x\in \partial K$, let $e_1(x),\ldots,e_N(x)$ be the unit vectors representing the orthonormal base of the coordinate system for which the previous representations hold. Then we assume that $e_N(x)$ is a Lipschitz function of $x\in\partial K$, with Lipschitz constant bounded by $L$, and $e_{N-1}(x)$ is a Lipschitz function of $x$, as $x$ varies in the boundary of $\partial K$, with Lipschitz constant bounded by $L$.

We remark that this implies that any connected component of $\partial K$ is an oriented Lipschitz hypersurface with or without boundary. In the case of a hypersurface without boundary, then this is the boundary of a Lipschitz domain  contained in $K$.

Fixed $h>0$, let
\begin{multline}\label{Khdef}
K_h=K\cup \{x+te_N(x):\ x\in \partial K,\ t\in[-h,h]\}\cup\\ \cup\{x+se_{N-1}(x)+te_N(x):\ x\text{ in the boundary of }\partial K,\ t\in[-h,h],\ s\in[0,h] \}.
\end{multline}

Notice that $K_h\subset B_{\tilde{h}}(K)$ for any $\tilde{h}>\sqrt{2}h$. In order to construct the function $\tilde{d}$ and 
prove the required properties, without loss of generality we may assume that $K$ has only one connected component.

We begin with the simpler case of $K=\overline{D}$, $D$ being a Lipschitz domain. In such a case,
at least for $h\leq r/2$,
$$
K_h=K\cup \{x+te_N(x):\ x\in \partial K,\ t\in[0,h]\}.
$$
By a standard use of the contraction mapping principle, it is not difficult to prove that there exist positive constants $h_0$, $c$, $0<c<1$, $r_1$ and $L_1$, depending on $r$ and $L$ only, such that for any $h$, $0<h\leq h_0$, the following holds.

For 
any $x\in\partial K_h$ there exists a function $\varphi^h:\mathbb{R}^{N-1}\to\mathbb{R}$, such that $\varphi(0)=0$ and which is Lipschitz with Lipschitz constant bounded by $L_1$, such that, up to a rigid change of coordinates, we have $x=0$ and
$$B_{r_1}(x)\cap K_h= \{y\in B_{r_1}(x):\ y_N\leq \varphi^h(y')\}.$$
Furthermore, for any $x\in\partial K$ we have that $B_{ch}(x)\subset K_h$, that is $B_{ch}(K)\subset K_h$. Finally $K_h$ and its boundary are connected and its exterior is connected as well.

We notice that such a Lipschitz condition implies that $B_{R+1}\backslash K_h$ satisfies a cone condition with a cone not depending on $h$ but only on $r$ and $L$. Therefore, by Lemma~\ref{conelemma}, the required properties are all satisfied.

A sketch of the proof is the following. Let us fix $x\in \partial K$. Without loss of generality, we take $x=0$ and assume that
$$K\cap B_r=\{y\in B_r:\ y_N\leq \varphi(y')\},$$
with $\varphi:\mathbb{R}^{N-1}\to\mathbb{R}$ a Lipschitz function with Lipschitz constant bounded by $L$ and such that
$\varphi(0)=0$. There exists a constant $c_1$, $0<c_1\leq 1/4$, depending on $L$ only such that the graph of $\varphi$
with $y'\in B'_{c_1r}$ is contained in $B_{r/4}$ and, consequently, it is a connected set.

For any $y'\in B'_{c_1r}$, let $\nu(y')=e_N(y',\varphi(y'))$, which is a Lipschitz function of $y'$ with Lipschitz constant bounded by $L(L+1)$. Notice that $\nu(0)=e_N(0)=e_N$ and that, as usual,
$\nu=(\nu',\nu_N)\in\mathbb{R}^{N-1}\times\mathbb{R}$.
Fixed $h$, let
$G_h:B'_{c_1r}\to \mathbb{R}^{N-1}$ be defined as follows. For any $y'\in B'_{c_1r}$, $G_h(y')=y'+h\nu'(y')$. Provided $|h|L(L+1)\leq 1/2$,
we have that $G_h$ is Lipschitz with Lipschitz constant bounded by $3/2$. Moreover, $G_h$ is injective and, by the contraction mapping principle, we have that $G_h(B'_{c_1r})$ contains $B'_{c_1r/2}$. Finally, $G_h^{-1}$ is Lipschitz on $B'_{c_1r/2}$ with a Lipschitz constant bounded by $2$.
For any $z'\in B'_{c_1r/2}$, let $\varphi^h(z')=h\nu_N(G_h^{-1}(z'))+\varphi(G_h^{-1}(z'))$. Such a function $\varphi^h$ is Lipschitz with a Lipschitz constant bounded by $L_1$, $L_1$ depending on $L$ only. Again we notice that $\varphi^h(0)=h$.

Next, we can find constants $h_0$, $0<h_0\leq 1/(2L(L+1))$, and $r_1$, $0<r_1\leq c_1r/2$,
depending on $r$ and $L$ only, such that for any $z'\in B'_{r_1}$
we have that $\varphi^h(z')$ is a strictly increasing function with respect to $h$, $-h_0\leq h\leq h_0$. Therefore, possibly after reducing $r_1$, for any $h>0$ and any $z'\in B'_{r_1}$
we have that
$$B_{r_1}(z',\varphi^h(z'))\cap \partial K_h=\{y\in B_{r_1}(z',\varphi^h(z')):\ y_N=\varphi^h(y')\}.$$
Thus we have obtained that $\partial K_h$ satisfies a Lipschitz condition, with constants depending on $r$ and $L$ only, and hence $\mathbb{R}^N\backslash K_h$ and $K_h$ satisfy a cone condition with a cone depending on $r$ and $L$ only.
Finally, taking $z'=0$, it is easy to show that there exists $c$, $0<c<1$ depending on such a cone only, such that $B_{ch}(0)\subset K_{h}$, that is $B_{ch}(K)\subset K_h$.

Then, for any $x\in K_{h_0}$, we define $\tilde{d}(x)=\min\{h\geq 0:\ x\in K_h\}$, where $K_0=K$.
For any $x\in\mathbb{R}^N\backslash K_{h_0}$, we define $\tilde{d}(x)=h_0+\mathrm{dist}(x,K_{h_0})$. Clearly the definition of $K_h$ with respect to the distance $\tilde{d}$ is consistent with
\eqref{Khdef} for any $h$, $0<h\leq h_0$. Finally, it is now easy to show that the function $\tilde{d}$ satisfies all the required properties.

In the case when $K$ is equal to $\partial K$ and consists of one connected component given by a Lipschitz hypersurface with boundary, we can prove that $B_{R+1}\backslash K_h$ satisfies a cone condition with a cone not depending on $h$, for any $h$ small enough. The arguments
and the construction of the function $\tilde{d}$
are similar to the one exposed above with a little extra care along the boundary of the hypersurface. We leave the lengthy but straightforward details to the reader.

\end{document}